# Frequency Domain Design of Fractional Order PID Controller for AVR System Using Chaotic Multi-objective Optimization

Indranil Pan[a], Saptarshi Das[b,c*]

a) Centre for Energy Studies, Indian Institute of Technology Delhi, Hauz Khas, New Delhi 110 016, India
b) Department of Power Engineering, Jadavpur University, Salt Lake Campus, LB-8, Sector 3, Kolkata-700098, India.
c) Communications, Signal Processing and Control Group, School of Electronics and Computer Science, University of Southampton, Southampton SO17 1BJ, United Kingdom.

Emails: indranil.jj@student.iitd.ac.in (I. Pan),
saptarshi@pe.jusl.ac.in, s.das@soton.ac.uk (S. Das*)

**Abstract**

A fractional order (FO) PID or FOPID controller is designed for an Automatic Voltage Regulator (AVR) system with the consideration of contradictory performance objectives. An improved evolutionary Non-dominated Sorting Genetic Algorithm (NSGA-II), augmented with a chaotic Henon map is used for the multi-objective optimization based design procedure. The Henon map as the random number generator outperforms the original NSGA-II algorithm and its Logistic map assisted version for obtaining a better design trade-off with an FOPID controller. The Pareto fronts showing the trade-offs between the different design objectives have also been shown for both the FOPID controller and the conventional PID controller to enunciate the relative merits and demerits of each. The design is done in frequency domain and hence stability and robustness of the design is automatically guaranteed unlike the other time domain optimization based controller design methods.

## 1. Introduction

Large power distribution networks must keep the overall voltage profiles at an acceptable level at all times. The connected equipments are designed for a particular nominal voltage and frequency of operation and any aberration from the nominal case generally leads to a decrease in performance and reduction in life time of these equipments. Frequent fluctuations in the load of the power network affects the voltage profile and hence the power utility companies employ a wide range of devices like capacitor banks, on-load tap changing transformers, automatic voltage regulators (AVRs) etc. [1–3] to keep the operational voltage profile at an acceptable level. Additionally, the amount of line losses due to the flow of real power depends on the reactive power which in turn depends on the system voltage. Hence, control of the system voltage is a crucial aspect in the effective operation of the power system. To alleviate these issues to some extent, the AVR is connected to the power generating plants. The AVR system maintains the terminal voltage of the alternator in the generating station and also helps in suitable distribution of the reactive power amongst the parallel connected generators [4].

Traditionally the PID controller has been used in the AVR loop due to its simplicity and ease of implementation [5]. However, recently the fractional order PID (FOPID) controller have been used in the design of AVR systems and have been shown to outperform the PID in many cases [6], [7]. In Zamani *et al.* [8], the FOPID has been tuned for an AVR system using the Particle Swarm Optimisation (PSO) algorithm employing time domain criterion like the Integral of Absolute Error (IAE), percentage overshoot, rise time, settling time, steady state error, controller effort etc. In Tang *et al.* [6], the optimal parameters of the FOPID controller for the AVR system, has been found using a chaotic ant swarm algorithm. In [6] a customised objective function has been designed using the peak overshoot, steady state error, rise time and the settling time. The above mentioned literatures perform optimisation considering only a single objective. But in a practical control system design multiple objectives need to be addressed. In the study by Pan and Das [9], the AVR design problem has been cast as a multi-objective problem and the efficacy of the PID and the FOPID controllers are compared with respect to different contradictory objective functions like the Integral of Time Multiplied Squared Error (ITSE) and the controller effort etc. However, the optimisation is done in the time domain and the obtained controller values are checked for robustness against gain variation by varying different parameters of the control loop. All these above mentioned literatures which employ time domain optimisation techniques cannot guarantee a certain degree of gain or phase margins which are important for the plant operator. These margins are useful from a control practitioner's view point as they can give an estimate of how much uncertainty the system can tolerate before it becomes unstable. Uncertainties can arise not only due to load variations in the power system, but there can be significant uncertainty due to modelling approximations or other stochastic phenomena. Hence frequency domain designs are mostly preferred over time domain design from the implementation and operation point of view of a control system. In spite of the importance of AVR in power systems, very few literatures consider a multi-objective formalism. A co-ordinated tuning of AVR and Power System Stabiliser (PSS) has been done in Viveros *et al.* [10] using the Strength Pareto Evolutionary Algorithm (SPEA). However, the contradictory objectives considered are the integrated time domain response for the AVR and the closed loop eigenvalue damping ratio of the PSS. This is a coupled time-frequency domain approach and does not address the inherent contradictory objectives in the AVR itself. In Ma *et al.* [11] a multi objective problem has been formulated for finding out the optimal solution for coordinate voltage control. A hierarchical genetic algorithm has been proposed for multi objective optimisation and a Pareto trade-off is obtained. In Mendoza *et al.* [12] a micro genetic algorithm is used to solve the multi objective problem of finding the AVR location in a radial distribution network in order to reduce energy losses and improve the energy quality. In [13], a similar problem has been attempted using a multi objective fuzzy adaptive PSO algorithm. However none of these papers consider the inherent design trade-off in the AVR tuning itself, which is one of the main focuses of the present paper.

In this paper an evolutionary multi-objective optimisation algorithm, the Non-dominated sorting genetic algorithm-II or NSGA-II [14], augmented with a chaotic Henon map, is used for designing a FOPID controller in frequency domain with contradictory objectives. The proposed frequency domain design methodology show that the FOPID controller is better than its PID counterpart for the considered set of objective functions. To the best of the author's knowledge, this is the first paper to make a comparative investigation into the multi-objective design trade-offs in frequency domain for the FOPID and the PID controller for an AVR system, using a chaotic map augmented multi-objective optimization algorithm.

The rest of the paper is organised as follows. Section 2 briefly introduces the concept of fractional calculus and the FOPID controller. In Section 3, the need for multi-objective optimisation, the description of the AVR system, the contradictory objective functions and the chaotic NSGA-II algorithm is discussed in detail. Section 4 illustrates the simulation results along with a few discussions. The paper ends in Section 5 with the conclusions followed by the references.

## 2. Fractional calculus and the fractional order PID (FOPID) controller

Fractional calculus is an extension of the integer order differentiation and integration for any arbitrary number. The fundamental operator representing the non-integer order differentiation and integration is given by ${}_aD_t^\alpha$ where $\alpha \in \mathbb{R}$ is the order of the differentiation or integration and $a$ and $t$ are the bounds of the operation. It is defined as

$${}_aD_t^\alpha = \begin{cases} \dfrac{d^\alpha}{dt^\alpha}, & \alpha > 0 \\ 1, & \alpha = 0 \\ \int_a^t (d\tau)^\alpha, & \alpha < 0 \end{cases} \quad (1)$$

There are three main definitions of fractional calculus, the Grünwald-Letnikov (GL), Riemann-Liouville (RL) and Caputo definitions. Other definitions like that of Weyl, Fourier, Cauchy, Abel and Nishimoto also exist. In the fractional order systems and control related literatures mostly the Caputo's fractional differentiation formula is referred. This typical definition of fractional derivative is generally used to derive fractional order transfer function models from fractional order ordinary differential equations with zero initial conditions. According to Caputo's definition, the $\alpha^{th}$ order derivative of a function $f(t)$ with respect to time is given by (2) and its Laplace transform can be represented as (3).

$$D^\alpha f(t) = \frac{1}{\Gamma(m-\alpha)} \int_0^t \frac{D^m f(t)}{(t-\tau)^{\alpha+1-m}} d\tau, \quad \alpha \in \mathbb{R}^+, m \in \mathbb{Z}^+ \quad (2)$$
$$m-1 \leq \alpha < m$$

$$\int_0^\infty e^{-st} D^\alpha f(t) dt = s^\alpha F(s) - \sum_{k=0}^{m-1} s^{\alpha-k-1} D^k f(0) \quad (3)$$

where, $\Gamma(\alpha) = \int_0^t e^{-t} t^{\alpha-1} dt$ is the Gamma function and $F(s) := \int_0^\infty e^{-st} f(t) dt$ is the Laplace transform of $f(t)$. This definition is used in the present paper for realizing the fractional integro-differential operators of the FOPID controller.

The fractional order PID controller is a generalization of its integer order counterpart where the integro-differential orders are two additional tuning knobs [15]. Thus in addition to the conventional proportional, integral and derivative gains $\{K_p, K_i, K_d\}$, there are also the integration and the differentiation orders $\{\lambda, \mu\}$. In the present study, $5^{th}$ order Oustaloup's recursive approximation is done for the integro-differential operators within a frequency band of the constant phase elements (CPEs) as $\omega \in \{10^{-2}, 10^{2}\}$ rad/sec. This frequency domain rational approximation method of realization is preferred over the others like Grunwald-Letnikov method since the realized approximate transfer functions can be easily implemented in real hardware using higher order Infinite Impulse Response (IIR) type analog or digital filters. The transfer function for the fractional order controller is given by Equation (4).

$$C(s) = K_p + \frac{K_i}{s^\lambda} + K_d s^\mu \tag{4}$$

## 3. Multi-objective optimization framework for FOPID controller design
### 3.1. Requirement for frequency domain multi-objective controller design

There are many controller design procedures like the $H_2$, $H_\infty$ or $L_1$ norm based designs where the controller design problem is reduced to that of minimizing the weighted norm of a closed loop transfer function. However, each of these norms addresses a specific performance criterion of the control system. For example in Herreros *et al.* [16], minimizing the $H_2$ norm implies closed loop stabilization in the presence of disturbances and minimizing $H_\infty$ norm gives closed loop robust stability. However, the control system designed with the $H_2$ norm minimization technique would have an arbitrary robustness as it has not been explicitly taken into the design criteria. Similarly for the $H_\infty$ norm case, the stabilization in the presence of disturbance is not addressed. In a practical control system design problem, the designer should design a system which ideally should have both properties to some extent. Hence the design algorithm must be capable of handling multiple objectives at the design phase itself. The NSGA-II algorithm is an evolutionary multi-objective optimisation algorithm which is suitable for designing such controllers with multiple objectives as shown in [9], [17], [18]. In [9], the multi objective design for fractional order controllers have been done by considering different conflicting time domain design criteria. However, the time domain design methods cannot explicitly quantify the parametric robustness of the designed system. The frequency domain design tools are much more powerful for analysing and designing linear systems. Hence, it is naturally advantageous to assimilate the different frequency domain design criteria in one design framework to study the various trade-offs of the design problem and choose the final controller based on these trade-offs.

### 3.2. Description of the AVR system

The salient components of the AVR system are the generator, exciter, amplifier and the sensor. Figure 1 shows the schematic diagram of the AVR with the fractional order PID controller. The output voltage of the generator $y(t)$, is repeatedly sensed by the voltage sensor and the signal is rectified, smoothed and analysed to see the deviation from the reference signal in the comparator. The error voltage resulting from this is used by the fractional order PID controller to generate a control signal which is then amplified and used to control the field windings of the generator by means of the exciter.

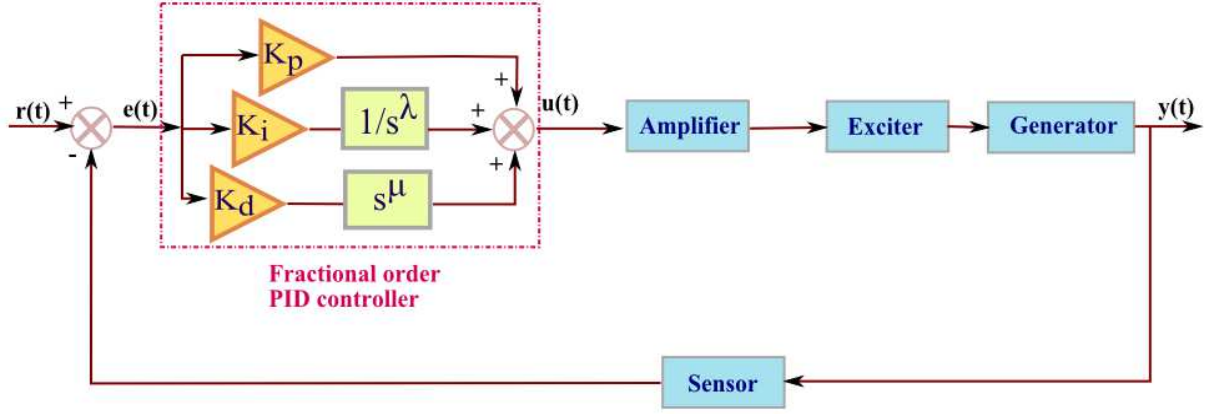

**Figure 1: Schematic of the AVR control loop with the fractional order PID controller**

The first order transfer function is able to capture the dynamics of the individual components of the AVR loop. This linear model takes care of the major time constants and neglects the saturation and other nonlinearities. The transfer functions of the different components of the system along with the range of their parameters are reported next [7]:

a) Amplifier model: $A(s) = \dfrac{K_A}{1+\tau_A s}$

   where $10 < K_A < 400$ and a small time constant is in the range $0.02 < \tau_A < 0.1$

b) Exciter model: $E(s) = \dfrac{K_E}{1+\tau_E s}$

   where $10 < K_E < 400$ and a time constant is in the range $0.5 < \tau_E < 1$

c) Generator model: $G(s) = \dfrac{K_G}{1+\tau_G s}$

   where $0.7 < K_G < 1$ and a time constant is in the range $1 < \tau_G < 2$. These constants vary depending on the load.

d) Sensor model: $S(s) = \dfrac{K_S}{1+\tau_S s}$

   where the time constant is in the range $0.001 < \tau_S < 0.06$

The values chosen for this are similar to those in [7]. Thus for the amplifier model $K_A = 10$ and $\tau_A = 0.1$. For the Exciter model $K_E = 1$ and $\tau_E = 0.4$. For the generator model $K_G = 1$ and $\tau_G = 1$. For the sensor model $K_S = 1$ and $\tau_S = 0.01$.

The inclusion of nonlinearities is definitely an important avenue towards real world modelling of the AVR system. Since this is a frequency domain design, the added nonlinearities make the system complicated and Bode analysis becomes cumbersome. This is more so with the fractional order controller than the integer order one. For the fractional order case, the mathematical formulation needs to be developed, to find out the stability region and the associated performances related to the AVR, since such studies do not exist till date. The objective of this paper is to direct the research towards deriving analytically tractable results for the fractional order case. This paper can be seen as a first step in establishing that a fractional order controller design in frequency domain can definitely give improvements among contradictory performance specifications over the conventional PID controller. The linearity assumptions as made in this paper is done in several others like [7], [19] and more recently in [8], [20–22].

Now, the closed loop transfer function ($G_{cl}(s)$) of the AVR system can be written in terms of the individual components as

$$G_{cl}(s) = \frac{C(s)A(s)E(s)G(s)}{1+C(s)A(s)E(s)G(s)} \tag{5}$$

Since the sensor is in the feedback loop, the frequency domain stability analysis tools are difficult to be applied for such a case. Hence for the closed loop transfer function in Equation (5), the effective open loop transfer function or $G_{ol\_eff}(s)$ is found out which is equivalent to Equation (5) with a unity feedback. Since both the transfer functions represent the same system, the following equations hold,

$$G_{cl}(s) = \frac{G_{ol\_eff}(s)}{1+G_{ol\_eff}(s)} \tag{6}$$

Thus the open loop effective transfer function with unity feedback can be calculated as

$$G_{ol\_eff}(s) = \frac{G_{cl}(s)}{1-G_{cl}(s)} \tag{7}$$

The bode plots and the gain and phase margins as reported in the next sections of the paper are all based on the effective open loop transfer function $G_{ol\_eff}(s)$ which is equivalent to the AVR system taken into consideration with a unity feedback instead of sensor transfer function being in the feedback path.

### *3.3. Conflicting objectives: Trade-off between gain crossover frequency and phase margin*

The gain crossover frequency ($\omega_{gc}$) and the phase margin ($\Phi_m$) are chosen as the two conflicting objectives for the design case, i.e.

$$\left. \begin{array}{l} J_1 = \omega_{gc} \\ J_2 = \Phi_m \end{array} \right\} \tag{8}$$

Here, the gain cross-over frequency and phase margins of the effective open loop system are related by the following equation

$$Arg\left[G_{ol\_eff}(j\omega)\right]_{\omega=\omega_{gc}} = -\pi + \Phi_m \qquad (9)$$

Both these objectives in (8) must be maximised for effective operation of the control loop. It is well known that high value of $\omega_{gc}$ makes the control system to act faster. Also with increase in speed the accuracy becomes low, which implies that a control system with high open loop gain and hence high $\omega_{gc}$ is prone to have oscillatory time response. These oscillations or overshoot is characterized by the phase margin and the relation between them is inversely proportional. Therefore it is logical that arbitrarily high speed ($\omega_{gc}$) and high accuracy in set point tracking ($\Phi_m$ or effective damping) cannot be obtained simultaneously for a linear control system. Increase in damping or $\Phi_m$ makes the system more sluggish and high speed of operation with increase in loop gain increases the overshoot. Therefore, it is imperative to study the trade-off between these two objectives as an effective measure of control system performance in frequency domain. Also, $\Phi_m$ is an important measure of robustness against system's gain variation which should be kept high. But too large $\Phi_m$ is undesirable as it results in a sluggish time response. For more details on fractional order controller design in frequency domain please look in [23][24].

Within the optimization algorithm while maximizing the conflicting objectives in (8), a constraint has been incorporated for the search with only those solutions yielding positive gain margin and phase margin; otherwise a large penalty is incorporated to discourage unstable solutions. It is well known that positive gain margin and phase margin of a linear system implies asymptotic stability (fractional order operators in the controller are linear as well). This is a necessary and sufficient condition for stability of the dominated and non-dominated solutions under nominal operating condition of the AVR system.

### *3.4. The Chaotic Non-dominated sorting Genetic Algorithm-II (NSGA-II)*

A generalized multi-objective optimization framework can be defined as follows [25], [26]:

$$\begin{aligned} &\text{Minimize } F(x) = (f_1(x), f_2(x),...,f_m(x)) \\ &\text{subject to: } g_i(x) \leq 0, \forall i \in [1, p], \\ &\qquad h_j(x) = 0, \forall j \in [1, q] \end{aligned} \qquad (10)$$

such that $x \in \Omega$

where $\Omega$ is the decision space, $\mathbb{R}^m$ is the objective space, $F : \Omega \to \mathbb{R}^m$ consists of $m$ real valued objective functions and $g_i(\cdot)$ and $h_j(\cdot)$ are the optional $p$ number of inequality and $q$ number of equality constraints on the problem respectively.

Let, $u = \{u_1,...,u_m\}, v = \{v_1,...,v_m\} \in \mathbb{R}^m$ be two vectors. Now, $u$ is said to dominate $v$ if $u_i < v_i \ \forall i \in \{1, 2,..., m\}$ and $u \neq v$. A point $x^* \in \Omega$ is called Pareto optimal if $\not\exists \ x \mid x \in \Omega$ such that $F(x)$ dominates $F(x^*)$. The set of all Pareto optimal points, denoted by PS is called the Pareto set. The set of all Pareto objective vectors, $PF = \{F(x) \in \mathbb{R}^m, x \in PS\}$, is called the Pareto Front. This implies that no other feasible objective vector exists which can improve one objective function without simultaneous worsening of some other objective function.

The NSGA-II algorithm [14] converts $m$ diverse objectives into one single fitness function by creating a number of different fronts. The solutions on these fronts are refined iteratively based on their distance with their neighbours (crowding distance) and their level of non-domination. The NSGA-II algorithm ensures that the solutions found are close to the original Pareto front and are diverse enough to find the whole length of the Pareto front.

Initially the algorithm starts with a population of randomly selected individuals from the search space. The individuals in the parent population of the NSGA-II algorithm are assigned a fitness value based on their non-domination level by checking the Pareto dominance. The non-dominated sorting algorithm is then used to allocate each solution on different fronts based on their domination level and the distance from the neighbouring solutions based on the crowding distance. The next generation is created from the parent generation using tournament selection and mutation. The pseudo code of the chaotic NSGA-II algorithm is outlined next in Algorithm 1.

**Algorithm 1**: The chaotic NSGA-II Algorithm

**Step 1. [Start]** Randomly generate N chromosomes within feasible search space

**Step 2. [Fitness]** Evaluate multiple fitness of each solution in the population

**Step 3. [Rank]** Rank Population by
    i. **[Domination rank]** Rank individuals using Algorithm 2
    ii. **[Crowding Distances]** Calculate using Algorithm 3

**Step 4. [New Population]** Create New Population through
    i. **[Selection]** Select parent chromosomes from previous population using crowding selection operator in Algorithm 4
    ii. **[Chaotic Crossover]** Using chaotic map and specified probability, do crossover of parents to form new offsprings.
    iii. **[Chaotic Mutation]** Using chaotic map and specified probability, mutate chromosomes to form new offsprings.
    iv. **[Accept]** Place new offsprings in the new population

**Step 5. [Replace]** Replace old population by new and continue

**Step 6. [Test]** If termination criteria is satisfied, return Pareto set of solutions from current population and stop

**Step 7. [Loop]** Goto **Step 2.**

The domination rank assignment algorithm as required in Algorithm 1 is outlined next.

**Algorithm 2**: Domination rank assignment

**Step 1.** Rank counter $r \leftarrow 0$
**Step 2.** $r \leftarrow r+1$
**Step 3.** Find non-dominated individuals from population $P$
**Step 4.** Assign rank $r$ to these individuals
**Step 5.** Remove these individuals from $P$ and continue
**Step 6.** If $P$ is empty then stop, else Goto **Step 2.**

The crowding distance represents the relative density of the solutions in the neighbourhood of a particular solution. Let a number of non-dominated solutions in $\Omega$ of size $\Psi$ be given along with a number of objective functions $f_k, k=1,2,\cdots,\Gamma$, where $\Gamma$ is the number of objectives. Let $d_i$ be the value of the crowding distance for solution $i$. Then $d_i$ can be calculated as given in Algorithm 3.

**Algorithm 3**: Crowding Distance assignment

**Step 1.** Let $d_i \leftarrow 0$ for $i=1,2,\cdots \Psi$
**Step 2.** For each objective function $f_k, k=1,2,\cdots,\Gamma$
  sort the set in ascending order
**Step 3.** Let $d_1$ and $d_\Psi$ be the maximum values, e.g. $d_1 = d_\Psi = \infty$
**Step 4.** For j=2 to $(\Psi-1)$, set $d_j \leftarrow d_j + \left(f_{k_{j+1}} - f_{k_{j-1}}\right)$

The crowding selection operator ($\succ$) is defined in Algorithm 4, which helps in comparing two chromosomes $x$ and $y$. Chromosome $x$ is better than $y$ if either of the two conditions given below is satisfied.

  i.   Domination rank of $x$ is smaller than $y$

  ii.  Domination ranks are equal and crowding distance of $x$ is larger than $y$

**Algorithm 4**: Crowding Selection Operator

$x \succ y$ iff $r_x < r_y$
  or $r_x == r_y$ AND $d_x > d_y$

In Caponetto *et al.* [27], an extensive study has been made to understand how different chaotic maps give better results for evolutionary algorithms. The study suggest that some discrete time chaotic sequences are able to give better results than random number generators, both in terms of objective function value and convergence speed. In Bucolo *et al.* [28], it has been documented that chaotic sequences give better optimisation results as compared to random number generators. This is because they introduce spatial diversity and non-organised patterns into the implementation of numerical procedures [28]. More recently, a lot of AVR literatures have focussed on intelligent single objective optimisation using chaotic maps for tuning PID controllers. In Zhu *et al.* [29], a chaotic ant swarm algorithm has been used to tune PID parameters and has been shown to give better results over the GA based tuning. Introducing other forms of randomness in the evolutionary or swarm algorithms have given better performance in tuning PID and PSS parameters in an AVR system [30]. Thus it has been extensively documented that chaotic sequences in single objective evolutionary optimisation give better results. In Guo *et al.* [31] it is shown that incorporating a chaotic process for the random number generation instead of the conventional random number generators, increases the efficiency of the algorithm and introduces diversity in the solutions. Hence in the present study, a Henon map is coupled with the standard NSGA-II algorithm to increase its effectiveness. The Henon map [32] is a two dimensional discrete time dynamical system that exhibits chaotic behaviour. Given a point with co-ordinates $\{x_n, y_n\}$, the Henon map transforms it to a new point $\{x_{n+1}, y_{n+1}\}$ using the following set of equations:

$$x_{n+1} = y_n + 1 - a x_n^2,$$
$$y_{n+1} = b x_n$$
(11)

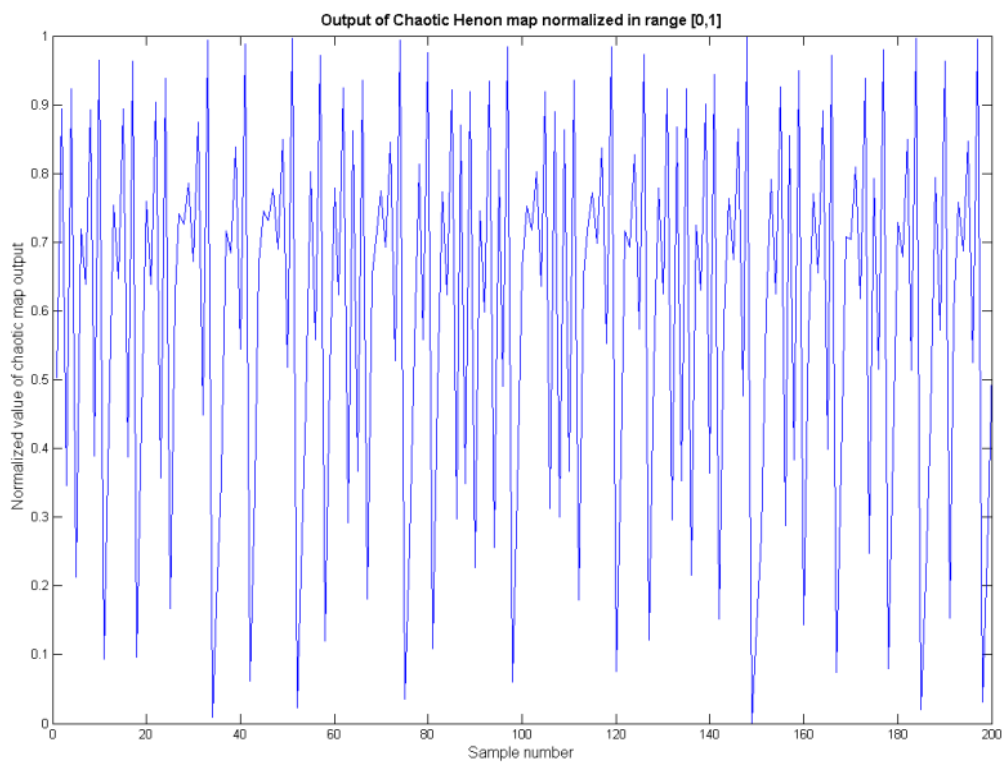

**Figure 2: Output of the Chaotic Henon map normalized in the range [0,1]**

The map is chaotic for the parameters $a = 1.4$ and $b = 0.3$. It is actually a simplified model of the Poincare section of the Lorenz system. The initial values of all the variables are zero. The output $y_{n+1}$ varies in the range $[-0.3854, 0.3819]$. Since the Henon map is used here as a random number generator, it must produce a random number in the range $[0,1]$. Hence the output is scaled in the range $[0,1]$ as also done in [33]. The first 200 samples of the scaled output of the Henon map are shown in Figure 2.

In order to get a fair comparison, in this paper we also compared the frequency domain multi-objective optimization results for FOPID controller tuning in AVR system with the classical NSGA-II and its chaotic versions augmented with Henon map and Logistic map. Pan and Das [9] studied time domain design trade-offs using Logistic map augmented NSGA-II in an AVR system. The one dimensional chaotic Logistic map is given as follows:

$$x_{n+1} = ax_{n+1}(1-x_n) \quad (12)$$

The initial condition of the map in Equation (12) has been chosen to be $x_0 = 0.2027$ and the parameter $a = 4$ has been taken similar to that in [27].

The population size for the chaotic Henon map augmented NSGA-II is taken as 200 and the number of generations as 150. The elite count, which represents the number of fittest individuals which are directly copied over to the next generation is taken as 30. An intermediate crossover scheme is adopted which produces off-springs by random weighted average of the parents. The mutation scheme adds a random number from a Gaussian distribution at an arbitrary point in the individual. For both the crossover and the mutation operations, the random numbers are generated from the chaotic Henon map among many other options like Logistic map etc. The search ranges of the controller parameters $\{K_p, K_i, K_d\}$ are chosen as $[0,10]$ and those of the fractional orders $\{\lambda, \mu\}$ are chosen as $[0,2]$. An in-house MATLAB code of the NSGA-II is developed for simulation purposes and coupled with the fractional order transfer function representations of the AVR system using the MATLAB scripting language.

Here the number of generations is relatively lesser compared to the population size. The population size is considered to be 200, so that more number of solutions can be obtained on the Pareto front and it is easier to visualize the different fronts and look at the appropriate trade-off. However, increasing the number of solutions also increases the computational cost and running the simulations for large number of generations is computationally intensive. This is more so with the case of Fractional order differ-integral operators, as they require the past history of the process in calculating the differ-integral value at any given point in time. Hence 150 generations are chosen, in which it is seen that the crowding distance between different individuals on the Pareto front does not change appreciably on increasing the number of generations. The crowding distance is essentially the Euclidean distance between different individuals on a front based on *m* different objectives in *m* dimensional hyperspace.

## 4. Results and discussions
### *4.1. Study of trade-off between frequency domain objectives*

Figure 3 shows the Pareto fronts for the FOPID and the PID controller showing the trade-offs between the gain crossover frequency ($\omega_{gc}$) and phase margin ($\Phi_m$) for maximization of both the objectives to achieve high speed and accuracy of control. As can be observed from the figure, the Pareto front for the PID totally lies in the region enclosed by the Pareto front of the FOPID controller and the X and Y axes. To understand the implication of this, let us take a constant phase margin line, of say 70 on the Y axis in Figure 3. The Pareto front for the PID and the FOPID controller intersect this horizontal line at approximately 8.4 radians/sec and 70.8 radians/sec respectively. Thus the FOPID controller is capable of giving a much higher gain crossover frequency than that with the PID. Consider a second case where the gain crossover frequency on the X axis is 10 rad/sec. For this frequency the phase margin of the PID controller on the Pareto front is around 60 and that of the FOPID controller is approximately 106. Thus the FOPID controller is capable of giving a much higher phase margin at the same gain crossover frequency. Also the Pareto front of the PID controller stops at 15 rads/sec whereas that of the FOPID controller continues over 110 rad/sec. Both the above mentioned examples are true for all the points on the Pareto front. Hence the FOPID controller outperforms the PID controller for frequency domain design of AVR system.

Compromise solution or other performance metrics can also be used to compare the obtained Pareto fronts. This is mostly done in cases where the improvement is small and cannot be captured visually, or if a benchmark needs to be set. In this case of controller design in AVR system however, it is easy to see the performance improvement visually from the Pareto front itself. The Pareto front of the FOPID controller completely dominates that of the PID controller in Figure 3, by a large margin. This implies that in all cases, the FOPID solutions would give better performance in both the objectives than the PID controller.

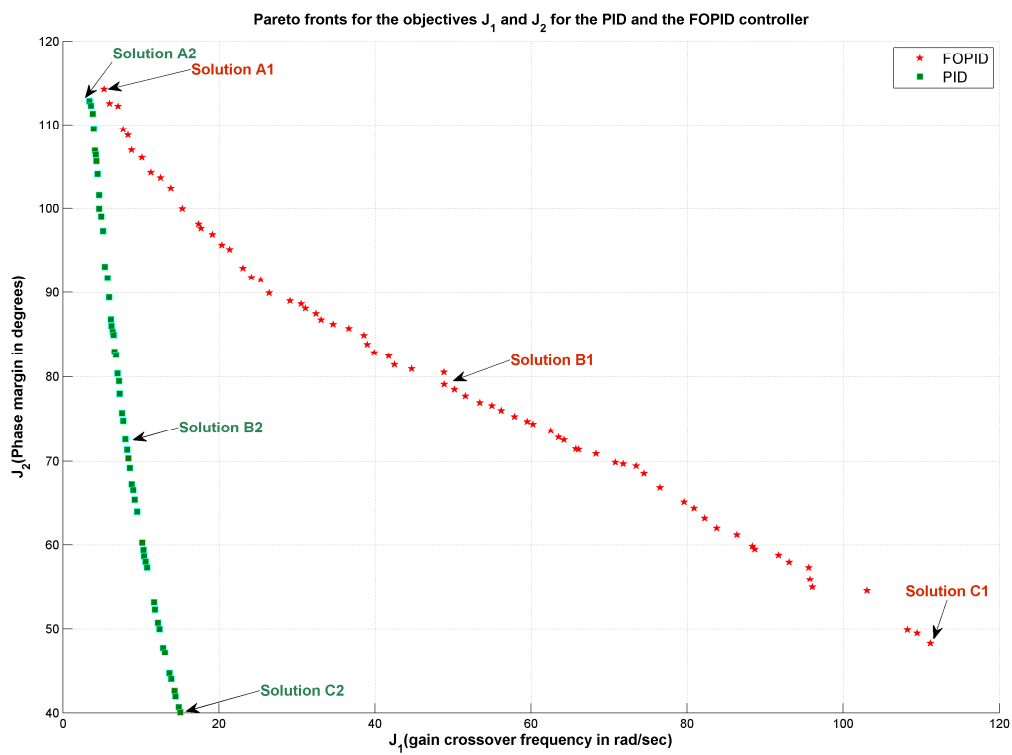

Figure 3: Comparison of the Pareto fronts for objectives $J_1$ and $J_2$ for the PID and the FOPID controllers

Table 1 shows some representative solutions as indicated on the Pareto front in Figure 3. The corresponding Bode diagrams for the PID and the FOPID controllers for these representative values are shown in Figure 4 and Figure 5 respectively. The trade-off between increasing in the gain crossover frequency leading to a decrease of the phase margin is clearly visible from these figures for both the PID and FOPID controller. But the FOPID maintains higher phase margin and hence higher damping and lower overshoot than that with the PID controller for similar increase in the speed of control that is characterised by the gain crossover frequency. A larger set of representative solutions on the Pareto front for the frequency domain trade-off between gain-crossover frequency and phase margin for PID and FOPID controllers has been given in the Appendix. Positive phase margins reported in the second column of Table 2 in the appendix confirms that all the Pareto optimal closed loop systems are stable for the nominal condition of the AVR (since the constraint imposed already guaranteed positive gain margin). The non-dominated solutions are also stable since the stability in terms of positive gain margin and phase margin have already been taken into consideration as a constraint before evaluation of the conflicting objective functions. The Pareto fraction for all simulation has been considered as 0.35. This represents the fraction of total solutions that lie on the Pareto front. Thus 70 non-dominated solutions are reported in the Appendix.

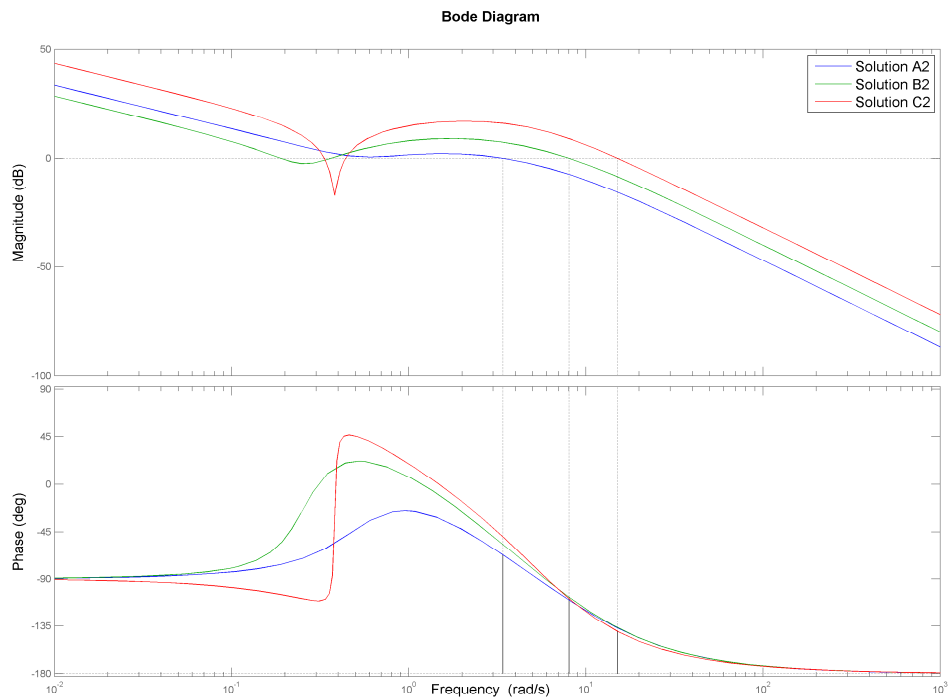

**Figure 4: Bode plots for the PID controller showing the gain crossover frequency and the phase margin**

A particular solution to the specific AVR problem is difficult to pin point without knowing the design constraints of the system. For example, in a particular case, there might be significant uncertainties in system identification of the AVR components. In such a case, the system designer might go for a more "safer" design, i.e. he can choose a solution which has more phase margin, so that even if the system parameters are different, the loop would not be unstable due to the safety factor considered by assigning a high value of phase margin to the system. However, as is evident, this would result in a consequent decrease in the performance of the system (reduction in fast damping of electro-mechanical oscillations). For a generic case, the median solutions $B_1/B_2$ may be chosen as an initial guess since it is somewhat balanced in both objectives, i.e., it has sufficient phase margin and appreciable gain crossover frequency.

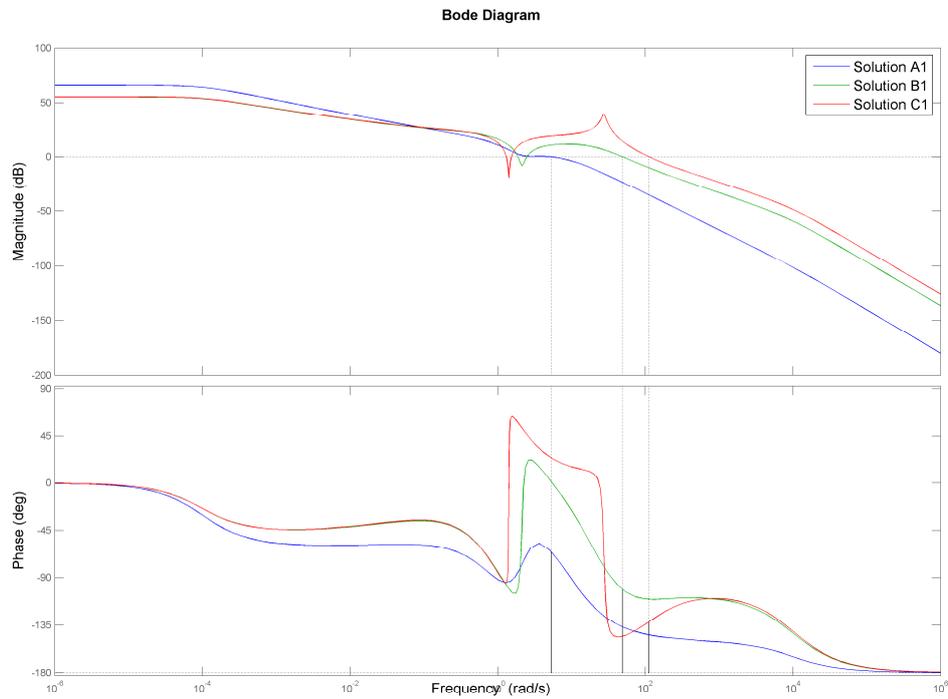

**Figure 5: Bode plots for the FOPID controller showing the gain crossover frequency and the phase margin**

**Table 1: Representative solutions on the Pareto front for the two controller types**

| Controller structure | Solution number | $J_1(\omega_{gc})$ | $J_2(\Phi_m)$ | $K_p$ | $K_i$ | $K_d$ | $\lambda$ | $\mu$ |
|---|---|---|---|---|---|---|---|---|
| FOPID | $A_1$ | 5.30366 | 114.23 | 0.408042 | 0.374094 | 0.17736 | 0.682778 | 1.333686 |
| | $B_1$ | 48.91055 | 79.09 | 0.963224 | 0.359946 | 0.281638 | 0.549116 | 1.830795 |
| | $C_1$ | 111.1908 | 48.28 | 1.037678 | 0.365733 | 0.654623 | 0.549738 | 1.871652 |
| PID | $A_2$ | 3.40643 | 112.82 | 0.124406 | 0.047737 | 0.180737 | - | - |
| | $B_2$ | 8.0412 | 72.58 | 0.077937 | 0.026848 | 0.40557 | - | - |
| | $C_2$ | 15.04809 | 40.07 | 0.015165 | 0.148481 | 1.018859 | - | - |

## *4.2. Robustness analysis of the obtained solutions*

For analysis of the robustness of the obtained solutions, two solutions are chosen from the PID and the FOPID Pareto fronts which have a phase margin of around 80 degrees as reported in Table 2. This typical choice makes both the PID and FOPID control loops to have a moderately high phase margin and low overshoot at nominal condition. Since the discrete solutions are taken from two different Pareto fronts, hence exactly the same value of phase margin is not obtained.

**Table 2: Controller values for testing robustness to gain variation**

| Controller type | Solution number | $J_1(\omega_{gc})$ | $J_2(\Phi_m)$ | $K_p$ | $K_i$ | $K_d$ | $\lambda$ | $\mu$ |
|---|---|---|---|---|---|---|---|---|
| FOPID | $D_1$ | 48.87964 | 80.48 | 1.010187 | 0.312872 | 0.268934 | 0.544849 | 1.841675 |
| PID | $D_2$ | 7.01888 | 80.34 | 0.059836 | 0.055538 | 0.345352 | - | - |

Now the robustness of the controllers is shown with respect to system gain variation. Figure 6 and Figure 7 show the effect of increase in the gain of the system. Since all the gains of the amplifier, exciter and the alternator are connected in series; the tolerance of gain for variation in one component would be the same if the variation were in the other components. Since the range of the exciter gains as given in Section 3.2 is much larger, we test the robustness of the tuned controllers as obtained in Table 2 by varying the gain of the exciter.

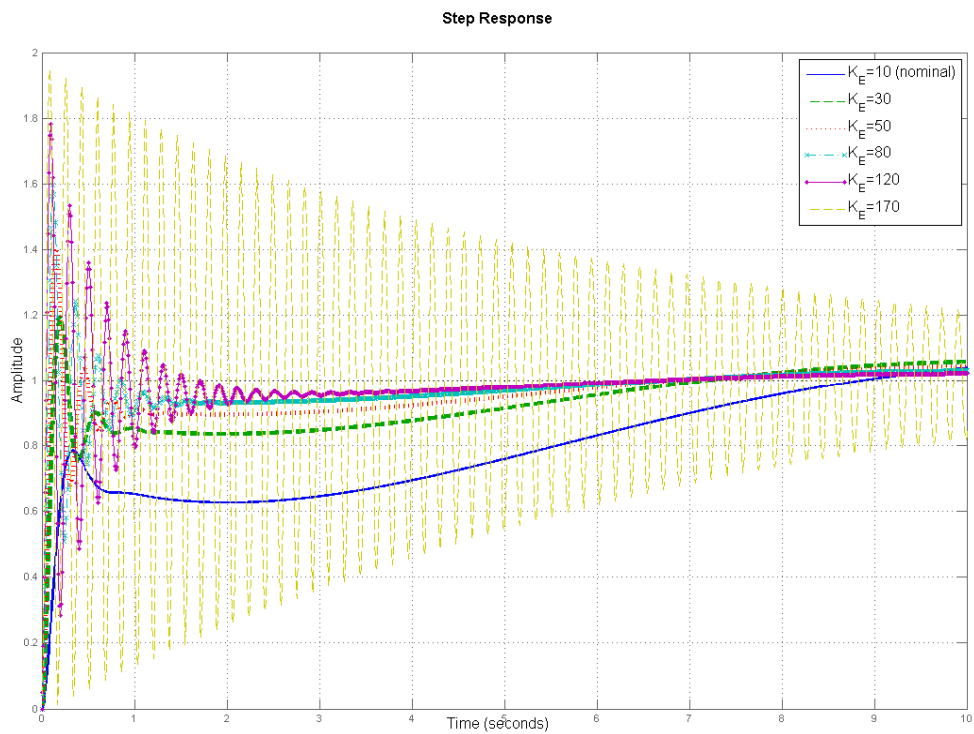

**Figure 6: Robustness of the PID controller for variation in system gain**

The time response of the PID controller for gain variation is shown in Figure 6. The nominal gain of the exciter is $K_E = 10$. The PID controller shows good response when $K_E$ is increased to 30 and 50. When $K_E$ is increased to 80, oscillations start to set in and for $K_E = 120$, high oscillations are observed. On increasing the gain further, the system shows unacceptably large oscillations and finally becomes unstable.

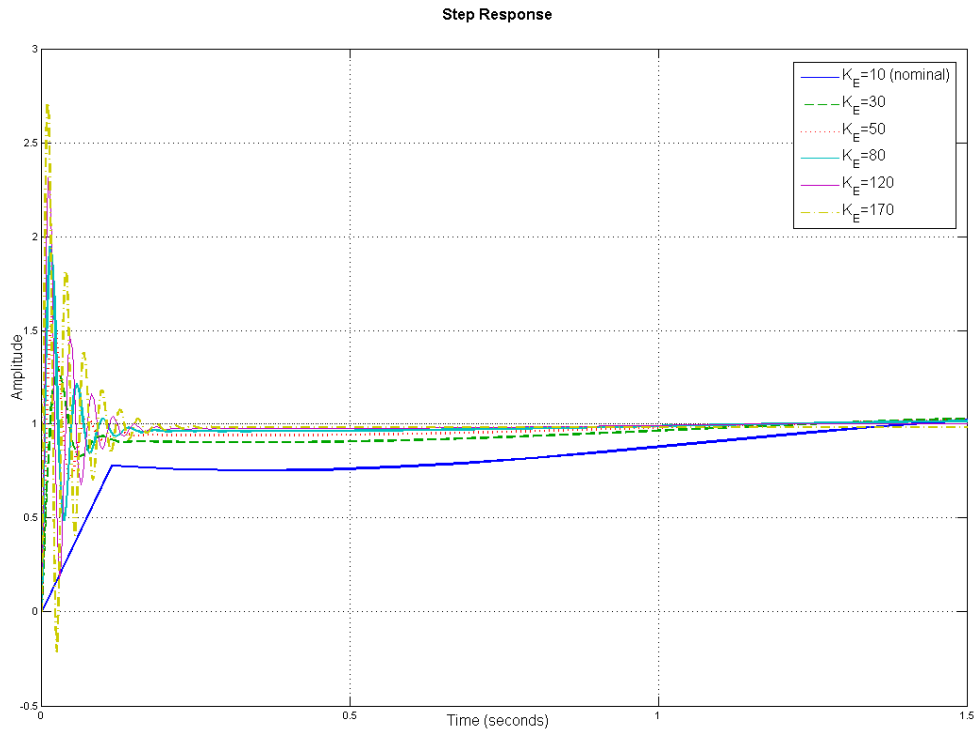

**Figure 7: Robustness of the FOPID controller for variation in system gain**

The time response of the FOPID controller is shown in Figure 7. It is clear that the FOPID gives a faster settling time than the PID controller. Since that order of the obtained FOPID controller is greater than 1, there is a derivative kick, at $t = 0$, when the set point is applied to the system. This derivative kick induces oscillations which die down very quickly (in 1/10th of a second). As the gain of the system is increased, the oscillations due to the effect of derivative kick increase in magnitude, but still die down in a very short time. This derivative kick can easily be removed while implementing it in practical hardware by incorporating the derivative action on the system output instead of the error which is a common practice [34]. The gain of the FOPID controller is increased to 170, and it can be seen that still the FOPID has good time response apart from the derivative kick effect whereas the PID controller shows very high oscillatory response in this case. Thus one important achievement of the present work is that, even if both the controllers are designed for almost same phase margin, the FOPID controller is capable of tolerating more gain variations than the PID controller, as has been shown in other fractional order controller design related literatures as well [24].

### 4.3. Comparison of design trade-offs with different chaotic map augmented NSGA-II algorithm and few discussions

As reported in [28], depending on the nature of the problem, different chaotic random number generators e.g. Logistic, Sinusoidal, Tent, Gauss, Lozi, Henon map, Lorenz and Chua system etc. can have different percentage improvements over the original algorithm. Since every problem setting is different, it is difficult to identify which chaotic map would give the best performance improvement. Comparing the performance improvement for all possible chaotic maps can be done through an exhaustive simulation study but that would be a digression from the main focus of the paper. But for the sake of completeness, performance comparison has been given in Figure 8 for the Logistic map and the Henon map augmented version of the NSGA-II algorithm with its classical version for studying the design trade-off for FOPID controller.

It is clear from Figure 8 that the original NSGA-II algorithm performs worse than its two chaotic versions. With the inclusion of Logistic map (12) for random number generation improves the results in terms of spread of the Pareto front as well as maximum achievable value of the conflicting objectives. The best result for the FOPID controller has been found using the Henon map (11) with a widely spread Pareto front and high value of phase margin and gain crossover frequency denoting high robustness and speed of operation of the AVR control loop. The comparison of the Pareto fronts in Figure 8 justifies the necessity of creating the randomness of the crossover and mutation operations of the evolutionary algorithm using a chaotic map rather than using a normally distributed random variable.

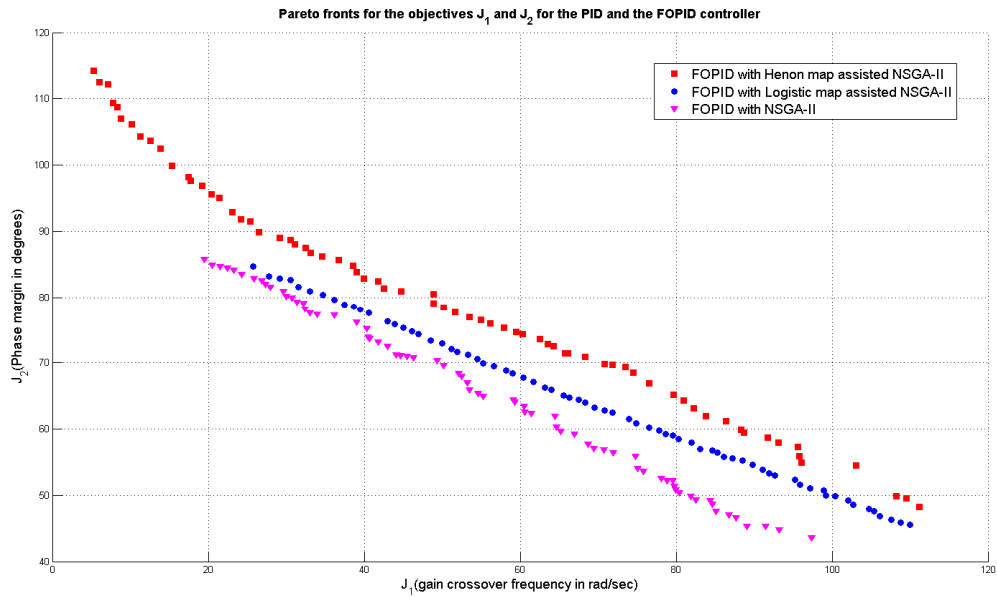

Figure 8: Comparison of the Pareto fronts for the FOPID controller with NSGA-II and its two chaotic versions

It may be argued that such a difference in the Pareto optimal solutions could also have been found by varying the crossover and mutation operators. Since, the main focus of the present paper is to compare the original NSGA-II and chaotic versions of NSGA-II for frequency domain FOPID controller design, the same mutation and crossover parameters have been specified for both algorithms for a fair comparison in Figure 8. Although for the sake of completeness and highlighting the effectiveness of chaos based random number generators in the mutation and crossover functions, the Pareto optimal fronts have been compared in Figure 9 for the chaotic Henon map augmented NSGA-II algorithm and its classical version. This comparison has been done with six different crossover fractions ($Cr$), varying between one and zero in a linear step of 0.2. In other words, this means the Mutation fraction being $M = 1 - Cr$, in each simulation. The Crossover and mutation operation are applied on the population in each generation excluding the elite members, which are directly carried on to the next generation.

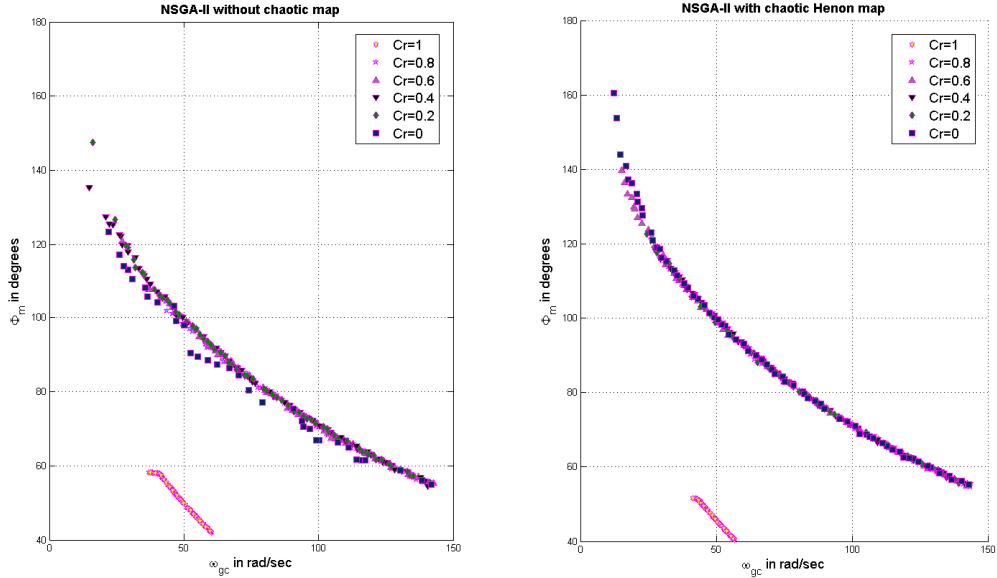

Figure 9: Effect of crossover fraction variation on the Pareto fronts of NSGA-II and its Henon map augmented version for FOPID controller in the AVR loop.

Experimentation with various values of mutation and crossover probabilities shows an interesting observation for replacing the random number generators with chaos. In Figure 9, it is observed that the chaotic Henon map augmented NSGA-II algorithm gives larger size of the Pareto front or a higher diversity of the population for $Cr = [0-0.8]$. For the extreme case of $Cr = 1$ when the whole population evolve through combination (crossover) of already obtained solutions in the first generation, both the algorithms give smaller and less diverse Pareto fronts, since the chance of exploring new solutions using the mutation function is not possible ($M = 0$). It is also observed from Figure 9 that variation in $Cr$ may produce larger Pareto front than that reported in Figure 8. But the improvement is mostly in the low gain crossover frequency ($\omega_{gc}$) and high phase margin ($\Phi_m$) regions which will result in slow operation of the AVR control loop and is not desirable from the practical point of view. We note that no improvement in the Pareto front is obtained where the solutions have both high speed (high $\omega_{gc}$) and low overshoot (high $\Phi_m$ and hence high damping). Also, the NSGA-II is an elitist algorithm and the number of individuals is taken as the same (30 elite solutions) for all the three versions of NSGA-II algorithm.

Therefore, it is observed that among different established techniques, one way of further improving the performance of powerful multi-objective evolutionary algorithms is to incorporate chaotic maps for the stochastic operations, instead of drawing the random numbers from a normal or uniform distribution, similar to the case of single objective cases in [28]. Since NSGA-II is a widely accepted multi-objective optimisation algorithm, therefore it has been chosen in this paper and has been improved by applying a chaotic map. The other newly introduced multi-objective optimization algorithms like Multi-objective Evolutionary Algorithm based on Decomposition (MOEA/D) has improvements in run time complexity over NSGA-II, but the performance improvements in obtaining better Pareto fronts is slight or almost similar [35]. The objective of the present paper is more applied in nature and geared towards the improvement in the power system control side using a harmonious blend of fractional order control techniques and multi-objective optimization based design trade-offs. The paper does not aim to propose a new algorithm or provide a comparison of a set of MOEA algorithms for test bench functions. Therefore, chaotic map augmented NSGA-II algorithm is highlighted to be used as a tool for an important frequency domain design aspect for AVR systems.

Also in recent literatures gain margin and phase margin based controller design have been attempted for various power system control related problem like hydro-power system [36], fuel-cell [37], harmonic compensation [38], maximum power point tracking [39], micro-grid [40], power factor correction [41], control of high voltage direct current (HVDC) links [42], sensitive loads [43] etc. In the above mentioned papers, the gain and phase margins are reported to represent robustness of the design though most of the controllers are designed in time domain. The present paper extends the concept of frequency domain design of fractional order controllers for an AVR system by maximizing both the conflicting objectives i.e. gain cross-over frequency and phase margin while also ensuring asymptotic stability using a constraint imposed as positive gain and phase margins.

The present paper advances the state-of-the-art techniques in both the fields of fractional order controller design in frequency domain and designing robust controller for AVR system. In a nutshell, the specific novel points of the paper over existing concepts are highlighted next:

- This is the first paper in which multi-objective design, in frequency domain controller tuning for the AVR system is taken up and an important control system design trade-off is shown.
- This is also the first paper to show a comparison of the design advantages of the FOPID controller over the PID controller in the AVR system, using a trade-off in frequency domain ($\omega_{gc}$ and $\Phi_m$ based).
- The coupling of MOEA for frequency domain designs for the fractional order PID controller is new as well.
- The chaotic Henon map augmented NSGA-II algorithm outperformed the classical version of the algorithm and its chaotic Logistic map augmented version for the design trade-off using FOPID controllers.
- Another important conclusion from the simulation results is that, even if both the controllers are designed at the same phase margin, the FOPID controller is capable of tolerating more gain variations than the PID controller. This insight in the AVR design could only be obtained due to frequency domain designs and not time domain ones as attempted previously by the contemporary researchers.

## 5. Conclusions

This paper investigates the design trade-offs between two frequency domain design specifications viz. phase margin and gain crossover frequency using a multi-objective formalism. The frequency domain design technique is more insightful than the time domain design methods since they give more knowledge about the stability and systems parametric robustness towards modelling uncertainties. A comparative analysis is made between the PID and the FOPID controller and it is shown that the latter outperforms the former and gives better designs. The inclusion of chaotic Henon map as random number generator for the mutation and crossover operations outperforms the original NSGA-II and its Logistic map assisted version for simultaneous maximization of phase margin and gain crossover frequency with a FOPID controller. It is also shown that even if both type of controllers are tuned at the same phase margin, the FOPID controller is capable of much faster time response or high gain crossover frequency while also tolerating more variations in system gain as compared to the PID controller. Future scope of work might be directed at multi-objective $H_2/H_\infty$ designs for such systems for better noise and disturbance rejection performance.

## Appendix

**Table 3: Additional representative solutions for the frequency domain multi-objective PID controller design**

| $J_1$ | $J_2$ | $K_p$ | $K_i$ | $K_d$ |
|---|---|---|---|---|
| 3.40643 | 112.82 | 0.124406 | 0.047737 | 0.180737 |
| 3.6007 | 112.28 | 0.104267 | 0.034093 | 0.187493 |
| 3.84237 | 111.31 | 0.078743 | 0.03586 | 0.197262 |
| 3.97081 | 109.52 | 0.08053 | 0.012235 | 0.200405 |
| 4.13343 | 106.92 | 0.088877 | 0.036938 | 0.207919 |
| 4.22759 | 106.46 | 0.078985 | 0.041031 | 0.212053 |
| 4.31403 | 105.66 | 0.075465 | 0.015707 | 0.214131 |
| 4.4711 | 104.16 | 0.069748 | 0.010661 | 0.220345 |
| 4.66527 | 101.65 | 0.074806 | 0.027171 | 0.22912 |
| 4.67043 | 99.91 | 0.104659 | 0.018351 | 0.228344 |
| 4.94505 | 99 | 0.066909 | 0.036714 | 0.241604 |
| 5.15112 | 97.31 | 0.056735 | 0.022062 | 0.250238 |
| 5.40485 | 93.06 | 0.09497 | 0.045057 | 0.262234 |
| 5.75968 | 91.73 | 0.043545 | 0.027328 | 0.278914 |
| 5.94787 | 89.41 | 0.058572 | 0.062203 | 0.289068 |
| 6.15745 | 86.78 | 0.081034 | 0.042326 | 0.298904 |
| 6.23253 | 86 | 0.085205 | 0.041784 | 0.302706 |
| 6.38986 | 85.3 | 0.065285 | 0.033469 | 0.310755 |
| 6.52623 | 84.95 | 0.038186 | 0.036623 | 0.318087 |
| 6.61373 | 82.91 | 0.083464 | 0.043393 | 0.322789 |
| 6.83377 | 82.56 | 0.031665 | 0.05081 | 0.335007 |
| 7.01888 | 80.34 | 0.059836 | 0.055538 | 0.345352 |
| 7.22379 | 79.47 | 0.030472 | 0.037499 | 0.356629 |
| 7.30076 | 77.99 | 0.068592 | 0.088957 | 0.362031 |
| 7.61399 | 75.69 | 0.069793 | 0.050558 | 0.379757 |
| 7.76649 | 74.82 | 0.059952 | 0.046848 | 0.388901 |
| 8.0412 | 72.58 | 0.077937 | 0.026848 | 0.40557 |
| 8.24005 | 71.32 | 0.075127 | 0.055222 | 0.418596 |
| 8.4271 | 70.32 | 0.062916 | 0.104176 | 0.431342 |
| 8.6109 | 69.17 | 0.061579 | 0.059699 | 0.442795 |
| 8.8213 | 67.27 | 0.100354 | 0.063242 | 0.457022 |
| 9.03266 | 66.58 | 0.062024 | 0.063659 | 0.471419 |
| 9.22402 | 65.45 | 0.061834 | 0.066189 | 0.484806 |
| 9.52636 | 63.91 | 0.046431 | 0.090341 | 0.50663 |
| 10.18044 | 60.26 | 0.050938 | 0.137505 | 0.555993 |
| 10.34469 | 59.41 | 0.050687 | 0.083392 | 0.56822 |
| 10.4264 | 58.65 | 0.082779 | 0.056202 | 0.574525 |
| 10.59847 | 58.03 | 0.058545 | 0.050422 | 0.588043 |
| 10.81551 | 57.34 | 0.01748 | 0.07462 | 0.605589 |
| 11.68309 | 53.12 | 0.027053 | 0.192735 | 0.680048 |
| 11.83177 | 52.27 | 0.0505 | 0.108613 | 0.692703 |
| 12.18815 | 50.7 | 0.052395 | 0.07083 | 0.724657 |
| 12.38919 | 49.95 | 0.036519 | 0.128848 | 0.743564 |
| 12.86585 | 47.73 | 0.077183 | 0.171238 | 0.789462 |
| 13.08361 | 47.14 | 0.031868 | 0.087765 | 0.809889 |
| 13.67029 | 44.67 | 0.070481 | 0.093846 | 0.869327 |
| 13.90453 | 44.01 | 0.03098 | 0.133941 | 0.893637 |
| 14.32036 | 42.57 | 0.018246 | 0.145279 | 0.937927 |
| 14.44865 | 41.94 | 0.055368 | 0.139746 | 0.952131 |
| 14.86378 | 40.67 | 0.019869 | 0.152144 | 0.998002 |
| 15.04809 | 40.07 | 0.015165 | 0.148481 | 1.018859 |
| 15.04809 | 40.07 | 0.015165 | 0.148481 | 1.018859 |

**Table 4: Additional representative solutions for the frequency domain multi-objective FOPID controller design**

| $J_1$ | $J_2$ | $K_p$ | $K_i$ | $K_d$ | $\lambda$ | $\mu$ |
|---|---|---|---|---|---|---|
| 5.30366 | 114.23 | 0.408042 | 0.374094 | 0.17736 | 0.682778 | 1.333686 |
| 6.00974 | 112.51 | 0.504434 | 0.361798 | 0.17661 | 0.641403 | 1.37913 |
| 7.11288 | 112.2 | 0.510327 | 0.35553 | 0.175536 | 0.614996 | 1.425375 |
| 7.75225 | 109.42 | 0.647101 | 0.374295 | 0.179113 | 0.67855 | 1.446619 |
| 8.36386 | 108.75 | 0.607594 | 0.41774 | 0.181423 | 0.536756 | 1.462657 |
| 8.8186 | 106.98 | 0.613146 | 0.424468 | 0.185917 | 0.525171 | 1.465211 |
| 10.15959 | 106.11 | 0.685626 | 0.455748 | 0.184696 | 0.512805 | 1.512731 |
| 11.30054 | 104.34 | 0.679635 | 0.381581 | 0.187251 | 0.576834 | 1.534293 |
| 12.54465 | 103.69 | 0.624996 | 0.362303 | 0.186014 | 0.623516 | 1.565701 |
| 13.85911 | 102.45 | 0.753987 | 0.378374 | 0.186933 | 0.597303 | 1.597441 |
| 15.33216 | 99.91 | 0.806349 | 0.415093 | 0.194535 | 0.533406 | 1.613974 |
| 17.4243 | 98.12 | 0.831157 | 0.385922 | 0.197744 | 0.546866 | 1.646323 |
| 17.73429 | 97.61 | 0.885538 | 0.383752 | 0.199544 | 0.556723 | 1.649049 |
| 19.21421 | 96.87 | 0.885599 | 0.378978 | 0.199974 | 0.55826 | 1.67192 |
| 20.38633 | 95.62 | 0.886868 | 0.375879 | 0.204245 | 0.549996 | 1.682298 |
| 21.41144 | 95.09 | 0.885893 | 0.344681 | 0.205005 | 0.562106 | 1.695299 |
| 23.0989 | 92.89 | 0.872653 | 0.36734 | 0.214702 | 0.55081 | 1.702554 |
| 24.17059 | 91.78 | 0.84917 | 0.328256 | 0.21961 | 0.563198 | 1.70852 |
| 25.37894 | 91.5 | 0.90893 | 0.345799 | 0.21932 | 0.555992 | 1.723352 |
| 26.48177 | 89.91 | 0.901594 | 0.368515 | 0.227933 | 0.5548 | 1.723914 |
| 29.14798 | 88.99 | 0.935232 | 0.37902 | 0.230106 | 0.548385 | 1.748285 |
| 30.5613 | 88.64 | 0.91352 | 0.356164 | 0.230455 | 0.544952 | 1.76083 |
| 31.13023 | 88.08 | 0.932163 | 0.364009 | 0.233607 | 0.535911 | 1.762089 |
| 32.45091 | 87.46 | 1.056362 | 0.343856 | 0.236182 | 0.641301 | 1.770573 |
| 33.12488 | 86.72 | 1.048703 | 0.394585 | 0.240743 | 0.596317 | 1.770823 |
| 34.65181 | 86.18 | 0.882714 | 0.364649 | 0.242621 | 0.529469 | 1.780411 |
| 36.68926 | 85.69 | 0.947139 | 0.324541 | 0.243535 | 0.571076 | 1.794515 |
| 38.59985 | 84.85 | 0.952478 | 0.320906 | 0.247411 | 0.525375 | 1.803525 |
| 39.03633 | 83.8 | 0.933042 | 0.356328 | 0.255143 | 0.540882 | 1.798401 |
| 39.93925 | 82.81 | 0.937358 | 0.362023 | 0.262079 | 0.541347 | 1.797395 |
| 41.77541 | 82.41 | 0.980197 | 0.310706 | 0.262774 | 0.577489 | 1.808287 |
| 42.52286 | 81.37 | 0.920179 | 0.367113 | 0.270654 | 0.547585 | 1.805327 |
| 44.73159 | 80.89 | 0.950974 | 0.332221 | 0.271657 | 0.540998 | 1.817267 |
| 48.87964 | 80.48 | 1.010187 | 0.312872 | 0.268934 | 0.544849 | 1.841675 |
| 48.91055 | 79.09 | 0.963224 | 0.359946 | 0.281638 | 0.549116 | 1.830795 |
| 50.19855 | 78.48 | 0.96523 | 0.358846 | 0.28535 | 0.549018 | 1.834102 |
| 51.62418 | 77.69 | 0.959991 | 0.341717 | 0.290651 | 0.567635 | 1.836666 |
| 53.46737 | 76.91 | 0.994219 | 0.366418 | 0.295259 | 0.554058 | 1.84164 |
| 54.9999 | 76.54 | 0.990475 | 0.358769 | 0.296331 | 0.554113 | 1.8477 |
| 56.19832 | 76 | 0.975895 | 0.346624 | 0.299793 | 0.54853 | 1.850293 |
| 57.90919 | 75.27 | 0.989516 | 0.348503 | 0.304275 | 0.560709 | 1.854236 |
| 59.52975 | 74.7 | 0.961327 | 0.330656 | 0.307364 | 0.579707 | 1.858656 |
| 60.28328 | 74.39 | 0.983498 | 0.347627 | 0.309307 | 0.553374 | 1.860313 |
| 62.54317 | 73.61 | 0.973921 | 0.345019 | 0.313714 | 0.521106 | 1.866093 |
| 63.54179 | 72.84 | 1.003619 | 0.354049 | 0.320361 | 0.551994 | 1.865379 |
| 64.25322 | 72.52 | 1.012827 | 0.341964 | 0.3227 | 0.543738 | 1.866506 |
| 65.76118 | 71.46 | 1.017175 | 0.354474 | 0.331988 | 0.558569 | 1.866053 |
| 66.09355 | 71.38 | 1.018228 | 0.354852 | 0.332288 | 0.56734 | 1.867086 |
| 68.37068 | 70.87 | 1.055566 | 0.317713 | 0.333955 | 0.671601 | 1.874248 |
| 70.82049 | 69.84 | 1.005594 | 0.328576 | 0.341733 | 0.610377 | 1.877957 |
| 71.82754 | 69.68 | 1.067714 | 0.327194 | 0.341716 | 0.621406 | 1.881401 |
| 73.50915 | 69.43 | 1.076121 | 0.301313 | 0.341516 | 0.672214 | 1.887123 |
| 74.48853 | 68.54 | 1.066385 | 0.323938 | 0.350768 | 0.714733 | 1.884865 |
| 76.53762 | 66.85 | 1.060252 | 0.33723 | 0.368926 | 0.611538 | 1.881264 |
| 79.63935 | 65.17 | 1.046772 | 0.345634 | 0.385888 | 0.625904 | 1.882099 |
| 80.92524 | 64.31 | 1.070374 | 0.333147 | 0.39582 | 0.651813 | 1.881066 |

| 82.28211 | 63.12 | 1.040741 | 0.342917 | 0.410954 | 0.603401 | 1.877873 |
| 83.81983 | 61.96 | 1.066683 | 0.357641 | 0.42598 | 0.647168 | 1.875569 |
| 86.40558 | 61.17 | 1.032547 | 0.363613 | 0.433265 | 0.581284 | 1.879987 |
| 88.38527 | 59.84 | 1.043527 | 0.365871 | 0.45165 | 0.60663 | 1.877863 |
| 88.69665 | 59.46 | 1.049191 | 0.348303 | 0.457543 | 0.596432 | 1.876315 |
| 91.75599 | 58.75 | 1.051391 | 0.363293 | 0.463563 | 0.607264 | 1.882552 |
| 93.11559 | 57.96 | 1.049719 | 0.359551 | 0.475044 | 0.575179 | 1.881783 |
| 95.62516 | 57.31 | 1.057709 | 0.3442 | 0.481807 | 0.609705 | 1.88602 |
| 95.78492 | 55.9 | 1.044708 | 0.342715 | 0.509201 | 0.596727 | 1.876336 |
| 96.0917 | 54.93 | 1.054582 | 0.389439 | 0.528401 | 0.626308 | 1.870427 |
| 103.097 | 54.52 | 1.051081 | 0.30031 | 0.521808 | 0.674523 | 1.89102 |
| 108.2668 | 49.9 | 0.996947 | 0.340228 | 0.617894 | 0.511692 | 1.874384 |
| 109.5115 | 49.51 | 0.997054 | 0.339315 | 0.625464 | 0.507196 | 1.875352 |
| 111.1908 | 48.28 | 1.037678 | 0.365733 | 0.654623 | 0.549738 | 1.871652 |